\documentclass{amsart}
\usepackage{graphicx}
\usepackage{array}
\usepackage{epsfig}
\usepackage{epstopdf}
\usepackage{amsmath, amsfonts,amssymb,amsthm}
\usepackage{subfigure}
\theoremstyle{plain}
\newtheorem{lem}{Lemma}[section]
\newtheorem{thm}[lem]{Theorem}
\newtheorem{prop}[lem]{Proposition}
\newtheorem{cor}[lem]{Corollary}

\theoremstyle{definition}
\newtheorem{defn}{Definition}[section]

\newtheorem{exmp}{Example}[section]

\begin{document}

\title{Linear relations for a generalized Tutte polynomial}
	
	\author{Gary Gordon}
\address{Dept. of Mathematics\\
	Lafayette College\\
   	Easton, PA  18042-1781}
 \email{gordong@lafayette.edu}
 

\begin{abstract}  Brylawski proved  the coefficients of the Tutte polynomial of a matroid satisfy a set of linear relations. We extend these relations to a generalization of the Tutte polynomial that includes greedoids and antimatroids. This leads to families of new identities for antimatroids, including trees, posets, chordal graphs and finite point sets in $\mathbb{R}^n$. It also gives a ``new'' linear relation for matroids that is implied by Brylawski's identities. 
\end{abstract}

\maketitle
\section{Introduction} The Tutte polynomial was defined by Tutte \cite{tutte} in the 1940's as a way to simultaneously encode the chromatic polynomial and the flow polynomial of a graph (equivalently, the chromatic polynomial of the dual when the graph is planar). The polynomial was generalized to matroids by Crapo \cite{cr} and Brylawski \cite{bry1} in the 1960's and 1970's, when substantial progress was made connecting the polynomial to an impressive collection of counting problems in graphs and matroids. The  expository chapter \cite{bryox} remains an excellent introduction and resource for the Tutte polynomial for graphs and matroids. For an interesting account of how Tutte came to work on this polynomial, see \cite{tutte1}.

Tutte's original definition used {\it basis activities} that depend on an ordering of the elements of the ground set of the matroid.  For a matroid $M$, this formulation expresses the polynomial  as a sum over the bases of the matroid: $T(M;x,y)=\sum_{i,j \geq 0}b_{i,j}x^iy^j $, where $b_{i,j}$ is the number of bases of the matroid having {\it internal activity} $i$ and {\it external activity} $j$. While the definitions of internal and external activities do not concern us here, it is worth pointing out that the $b_{i,j}$ are matroid invariants, not depending on any ordering of the ground set. (An alternative approach to basis activities for the Tutte polynomial of a graph was undertaken by Bernardi \cite{bern}. His use of activities allows a different interpretation for the $b_{i,j}$, still dependent on an order of the edges of the graph, and does not seem to extend to matroids.)

Brylawski \cite{bry2} discovered that the coefficients of the Tutte polynomial $T(M;x,y)$ of a matroid $M$ satisfy the following collection of linear relations. (There are other boundary conditions these coefficients satisfy; Brylawski gives a basis for all such relations that we give as  Theorem~\ref{T:basis} here.) Our main theorem generalizes these relations to {\it ranked sets}, a natural generalization of matroids and greedoids based on the rank function.

\bigskip
{\bf Theorem 3.1}  {\it Let $G=(S,r)$ be a ranked set with $|S|=n$.  Let $S(G;u,v)=\sum_{A \subseteq E}u^{r(S)-r(A)}v^{|A|-r(A)}$ be the corank-nullity generating function, and let $$T(G;x,y)=S(G;x-1,y-1)$$ be the Tutte polynomial of $G$. Write $T(G;x,y)=\sum_{i, j \in \mathbb{Z}} b_{i,j}x^i y^j$. Then}
\begin{enumerate}
\item {\it For all} $0 \leq k < n$,
$$\sum_{i=0}^k \sum_{j=0}^{k-i} (-1)^j {k-i \choose j}b_{i,j}=0.$$
\item {\it For} $k=n$,
$$\sum_{i=0}^k \sum_{j=0}^{k-i} (-1)^j {k-i \choose j}b_{i,j}=(-1)^{n-r(S)}.$$
\end{enumerate}

We list the first four relations in Table~\ref{Ta:affine}, simplifying as much as possible. Note that the second relation listed ensures that the coefficient of $x$ equals the coefficient of $y$, provided the set has at least 2 elements. The coefficient $b_{1,0}$ is the {\it beta} invariant, and it could equally well be taken to be $b_{0,1}$. It is not difficult to prove that $b_{1,0}>0$ for a matroid if and only if $M$ is a connected matroid (assuming $M$ has at least 2 points), a result originally due to Crapo \cite{cr1}.
\begin{table}[htdp]
\begin{center}
\begin{tabular}{ccc} \hline
$b_{0,0}$ & $=$ & $0$ \\
$b_{0,1}$ & $=$ & $b_{1,0}$ \\
$b_{0,2}+b_{2,0}$& $=$ & $b_{1,0}+b_{1,1}$ \\
$b_{3,0}+ b_{0,2}+b_{1,2}$& $=$ & $b_{2,0}+b_{0,3}$ \\ \hline
   \end{tabular}
\end{center}
\medskip
\caption{Simplified versions of the first four affine relations.}
\label{Ta:affine}
\end{table}%

A deeper connection between the Tutte polynomial and matroid connectivity was proven in \cite{me1}. The authors prove that  a matroid is connected if and only if its Tutte polynomial is irreducible in the polynomial ring $\mathbb{Z}[x,y]$. The proof is based on an analysis of the identities of Theorem~\ref{T:basis}.

 Brylawski published two different proofs of these identities. In Theorem 6.6 of \cite{bry2}, he counts the number of flats in a matroid of a given corank and nullity. His later proof (Prop. 6.3.2.f of \cite{bry3}) was based on more general properties involving the relation between the Tutte polynomial and the {\it corank-nullity} generating function. This latter proof is still matroid-theoretic, but it is possible to remove the dependence of the proof on the properties of the matroid rank function. This proof is the foundation for our proof of Thm.~\ref{T:affine} here.

Our motivation in extending Brylawski's relations to non-matroidal combinatorial structures is three-fold.
\begin{enumerate}
\item The Tutte polynomial has important connections with combinatorial objects that do not have a matroidal structure. Extending the Tutte polynomial to greedoids and antimatroids is the focus of \cite{cg,gm1,gposet,gm2}. The polynomial and some  specializations (a one-variable characteristic polynomial and the beta invariant) have revealed interesting combinatorial properties of the (non-matroidal) object under consideration. For instance, the generalized Tutte polynomial of a rooted tree is essentially a generating function for the number of rooted subtrees with $k$ edges and $l$ leaves. Then Theorem 2.8 of \cite{gm1} shows that two rooted trees have the same Tutte polynomial if and only if they are isomorphic, i.e., the polynomial is a complete invariant.  Thus, a rooted tree can be uniquely reconstructed from its edge-leaf subtree data. Understanding the relations satisfied by the Tutte polynomial coefficients in this context should shed light on these combinatorial structures.
\item Applying the family of identities to classes of antimatroids should be especially fruitful. In this case, we can use an expansion of the polynomial  from \cite{gposet} in terms of convex sets. Then we can express Tutte polynomial coefficients as a sum of $a_{i,j}$, where $a_{i,j}$ is the number of convex sets of size $i$ with exactly $j$ ``interior'' elements. We interpret ``interior'' in a variety of ways, depending on the antimatroid under consideration. This allows us to develop new identities for a variety of combinatorial objects that do not form a matroid in a natural way, including  finite subsets of $\mathbb{R}^n$, posets (the {\it double shelling poset antimatroid}), trees,  and chordal graphs (the {\it simplicial shelling antimatroid}). One example of such an identity is the following:

\medskip
\noindent {\bf Corollary \ref{C:pointset}} Let $S$ be a finite subset of $\mathbb{R}^n$. Let $\mathcal{C}_1$ be the collection of all convex sets with exactly one interior point. Then
$$\sum_{C\in\mathcal{C}_1} (-1)^{|C|} =(-1)^{n}|int(S)|.$$

\medskip

\item In a larger sense, this work attempts to isolate the properties of the Tutte polynomial that depend on the underlying matroid structure. Some of the elementary properties the coefficients satisfy for matroids are not true in our more general setting. For instance, $b_{i,j}<0$ is possible for greedoids, something that is impossible for matroids (where the $b_{i,j}$ count the number of bases with certain activities). This is made explicit in Example~\ref{E:badiden}.
\end{enumerate}

Our approach  also leads to a ``new'' result for matroids, and what appears to be a new combinatorial identity. The matroid result is the identity for the case $k=n$, given in Cor.~\ref{C:matroid}:

\medskip
\noindent
{\bf Corollary~\ref{C:matroid}} Let $M$ be a rank $r$ matroid on $n$ points.  Then if $T(M)=\sum b_{i,j}x^iy^j$, 
$$\sum_{i=0}^n \sum_{j=0}^{n-i} (-1)^j {n-i \choose j}b_{i,j}=(-1)^{n-r}.$$
\medskip

This property of the coefficients of the Tutte polynomial follows from the relations found by Brylawski. Indeed, Brylawski shows that the relations of Theorem~\ref{T:basis} form a basis for all linear relations the coefficients satisfy, so the above relation must follow from those listed in Theorem~\ref{T:basis}. But this appears not to have been noticed before.

A byproduct of our approach is a combinatorial identity Cor.~\ref{C:iden}. 

\medskip
\noindent
{\bf Corollary~\ref{C:iden}}: For all $k, n \geq 0$, 
$$\sum_{i=0}^n \sum_{j=0}^{n-i} (-1)^i {n-i \choose j} {k\choose i}{k \choose j}=1.$$
\medskip

It is a straightforward exercise to prove this (using Vandermonde convolution, for example), but it may also be new. Our proof follows from one of the lemmas (Lemma~\ref{L:power}) we use to prove the main theorem.

The paper is structured as follows. We introduce the definitions and fundamental results we will need in Section~\ref{S:defs}. The main result of the paper, Theorem~\ref{T:affine}, is proven in Section~\ref{S:proof}. The proof follows from several  lemmas, but the logical dependence is straightforward. Section~\ref{S:anti} gives interpretations of these identities for antimatroids, concentrating specifically on trees, posets, chordal graphs and finite subsets of $\mathbb{R}^n$. Finally, we give examples and counterexamples in Section~\ref{S:ex}. Several examples demonstrate how much the identities given here determine the Tutte polynomial of a matroid or greedoid, but also show that the Tutte polynomial does not distinguish matroids from greedoids.

The author is indebted to Tom Brylawski (1944 -- 2007), who has influenced and inspired much of this work. A close read of \cite{bry3} continues to reveal ideas that deserve deeper exploration.

\section{Definitions and preliminary results}\label{S:defs}
Many of the definitions and preliminary results we will need appear in \cite{gdual}; see that reference for proofs. 

\begin{defn}\label{D:ranked set}
A {\it ranked set} is a set $S$ with a rank function $r$. We write $G=(S,r)$, where the function $r:S \to \mathbb{Z}$ satisfies
\begin{itemize}
\item  [(R0)] $r(\emptyset)=0$ [{\it normalization}]
\item [(R1)] $r(A) \leq r(S)$ for all $A \subseteq S$ [{\it rank $S$ maximum}]
\item [(R2)] $r(A) \leq |A|$ for all $A \subseteq S$ [{\it subcardinality}]
\end{itemize}
\end{defn}

We define the generalized Tutte polynomial $T(G;x,y)$ using the rank function:

\begin{defn} Let $G=(S,r)$ be a ranked set. Then the {\it generalized Tutte polynomial} is
$$T(G;x,y)=\sum_{A \subseteq S}(x-1)^{r(S)-r(A)}(y-1)^{|A|-r(A)|}.$$

\end{defn}
Properties (R1) and (R2) ensure the generating function we use to define $T(G;x,y)$ is, in fact, a polynomial. 
 Ranked sets generalize matroids and greedoids. Matroids and greedoids can be defined as ranked sets whose rank functions satisfy properties more restrictive than (R0), (R1) and (R2).

\begin{defn} A \textit{matroid} $M$ is a pair $(S,r)$ where $S$ is a finite set and $r:2^S \to \mathbb{Z}^+\cup\{0\}$ such that:
\begin{enumerate}
\item $r(\emptyset)=0$ [\textit{normalization}]
\item $r(A) \leq r(A \cup p) \leq r(A)+1$ [\textit{unit rank increase}]
\item $r(A\cap B)+r(A \cup B) \leq r(A)+r(B)$ [\textit{semimodularity}]
\end{enumerate}
\end{defn}
A standard reference for matroids is Oxley's text \cite{ox}. Greedoids were introduced by Korte and Lovasz \cite{kl} as a generalization of matroids in an attempt to isolate structures where the greedy algorithm always produces an optimal solution.  An extensive introduction appears in \cite{bz}.  Although there are fewer axiomatizations of greedoids than there are of matroids, it is also possible define greedoids from a  rank function.

\begin{defn}\label{D:greedoid} A {\it greedoid} $G$ is a pair $(S,r)$ where $S$ is a finite set and $r:2^S \to \mathbb{Z}^+\cup\{0\}$ such that:
\begin{enumerate}
\item $r(\emptyset)=0$ [{\it normalization}]
\item $r(A) \leq r(A \cup \{p\})$ [{\it increasing}]
\item $r(A) \leq |A|$ [{\it subcardinality}]
\item  If $r(A)=r(A \cup p_1)=r(A\cup p_2)$, then $r(A\cup\{p_1,p_2\})=r(A).$ [{\it local semimodularity}]
\end{enumerate}
\end{defn}

\begin{defn}\label{D:gens} Let $G=(S,r)$ be a ranked set.  Define duality, deletion and contraction:
\begin{itemize}
\item \textit{Duality} $G^*=(S,r^*)$, where $r^*(A):=|A|+r(S-A)-r(S)$.
\item \textit{Deletion} $G-p=(S-p,r')$, where $r'(A)=r(A)$ for all $A \subseteq S-p$.
\item \textit{Contraction} $G/p:=(G^*-p)^*$.
\end{itemize}
\end{defn}

These generalize the usual matroid definitions. Note that $r^*(S)=|S|-r(S)$.  This definition of duality also has the  involution property:  $(G^*)^*=G$.  
%

The next result is useful for computing the rank of a subset in $G/p$, and will also be needed in our proof of a deletion-contraction recursion for the Tutte polynomial (Theorem~\ref{T:tutte-gen}(1)).  For clarity, we may denote the rank function of $G$ by $r_G$.

\begin{thm}[Theorem 2.2 \cite{gdual}] \label{T:contract-rank}  Let $G=(S,r)$, where $r:2^S\to \mathbb{Z}$  satisfies $r(\emptyset)=0$.  Then, for all $p \in S$ and $A \subseteq S-p$, 
$$r_{G/p}(A)=r(A \cup p)-r(p).$$
\end{thm}

Applying Definition~\ref{D:gens} to the polynomial, we  get the following.

\begin{thm}[Theorem 3.1 \cite{gdual}] \label{T:tutte-gen} Let $r:2^S \to \mathbb{Z}$ be any function satisfying $r(\emptyset)=0$, and let $G=(S,r)$ be a ranked set.  Then 
\begin{enumerate}
  \item Deletion-contraction:  For any $p \in S$, $$f(G;t,z)=t^{r(G)-r(G-p)}f(G-p;t,z)+z^{1-r(p)}f(G/p;t,z).$$
  \item Duality:  Let $G^*$ be the dual of $G$ in the sense of Definition~\ref{D:gens}.  Then $$f(G^*;t,z)=f(G;z,t).$$ 
\end{enumerate}

\end{thm}

We will need a generalization of the matroid operation of {\it truncation}.

\begin{defn}\label{D:trunc}
Let $G=(S,r)$ be a ranked set.  Define the {\it truncation} $T(G)=(S,r_T)$ by specifying its rank function: 
$$ r_TA)= \left\{ \begin{array}{lll}
r(A)& \mbox{if}& r(A) \leq r(S)-1, \\
r(A)-1 &\mbox{if} & r(A)=r(S). \\
\end{array} \right.$$
\end{defn}

We will also need the following construction, generalizing the matroid operations of free  extension and co-extension.

\begin{defn}\label{D:exten}
Let $G=(S,r)$ be a ranked set.  Suppose $p\notin S$.
\begin{enumerate}
\item Free extension:  $G+p$ is defined on the pair $(S \cup p, r')$ where
$$ r'(A)= \left\{ \begin{array}{lll}
r(A)& & p \notin A, \\
r(A) & & p \in A \mbox{ and } r(A)=r(S), \\
r(A)+1  & & p \in A \mbox{ and } r(A)<r(S),
\end{array} \right.$$
\item Free co-extension:  $G\times p=(G^*+p)^*$.
\end{enumerate}
\end{defn}

\begin{prop}\label{P:exten}  Let $G=(S,r)$ be a ranked set.
\begin{enumerate}
\item Rank: 
\begin{enumerate}
\item $r(G+ p)=r(G)$ 
\item $r(G\times p)=r(G)+1$,
\item  $r((G +p)/p)=r(G)-1$,
\item $r((G\times p)-p)=r(G)+1$.
\end{enumerate}

\item Deletion in free extension:  $(G+p)-p=G$.
\item Contraction in free co-extension:   $(G \times p)/p=G$.
\end{enumerate}

\end{prop}
\proof

\begin{enumerate}
\item \begin{enumerate}
\item This is immediate from the definition.
\item Let $|S|=n$.  Recall $r(G^*)=n-r(G)$.  From Definition~\ref{D:gens}(1), we have 
\begin{eqnarray*}
r(G\times p)&=&r((G^*+p)^*) \\
&=&r^*(G^*+p)\\
&=& n+1-r(G^*) \\
&=&n+1-(n-r(G)) \\
&=& r(G)+1.
\end{eqnarray*}

\item  This follows from Theorem~\ref{T:contract-rank}.
\item Note that $((G\times p)-p)^*=(G\times p)^*/p=(G^*+p)/p$.  By part (c), we know $r((G^*+p)/p)=r(G^*)-1$.  Hence
\begin{eqnarray*}
r((G\times p)-p) &=&n-r(((G\times p)-p)^*) \\
&=&n-(r(G^*)-1) \\
&=&n-(n-r(G)-1)\\
&=&r(G)+1.
\end{eqnarray*}
\end{enumerate}
\item  This is immediate from the definition.
\item This follows from part (2) and duality.
\end{enumerate}

\qed

We remark that the proof of Prop.~\ref{P:exten} is the same as the standard proof for matroids. 

\section{Proof of the main theorem}\label{S:proof}
In this section we prove our main theorem.
\begin{thm}\label{T:affine}  Let $G=(S,r)$ be a ranked set with $|S|=n$.  Let $S(G;u,v)=\sum_{A \subseteq E}u^{r(S)-r(A)}v^{|A|-r(A)}$ be the corank-nullity generating function, and let $$T(G;x,y)=S(G;x-1,y-1)$$ be the Tutte polynomial of $G$. Write $T(G;x,y)=\sum_{i, j \in \mathbb{Z}} b_{i,j}x^i y^j$. Then
\begin{enumerate}
\item For all $0 \leq k < n$,
$$\sum_{i=0}^k \sum_{j=0}^{k-i} (-1)^j {k-i \choose j}b_{i,j}=0.$$
\item  For $k=n$,
$$\sum_{i=0}^k \sum_{j=0}^{k-i} (-1)^j {k-i \choose j}b_{i,j}=(-1)^{n-r(S)}.$$
\end{enumerate}\end{thm}

The proof of Theorem~\ref{T:affine} will follow several lemmas.  We begin by investigating these identities for arbitrary 2-variable polynomials under shifting operations (multiplication by $x$ or $y$).

\begin{lem}\label{L:shift}  Let $p(x,y)=\sum_{i,j \geq 0}a_{i,j}x^i y^j$ be a two-variable polynomial with  coefficients in a commutative ring and let $I_k(p(x,y))=\sum_{i=0}^k \sum_{j=0}^{k-i} (-1)^j {k-i \choose j}a_{i,j}.$  
\begin{enumerate}
\item $I_k(x^m p(x,y))=I_{k-m}(p(x,y))$ for $0 \leq m \leq k$.
\item $I_k(y \cdot p(x,y))=(-1)\sum_{s=0}^{k-1} I_{s}(p(x,y))$.
\end{enumerate}
\end{lem}

\proof  We prove (1) for $m=1$; the result follows by induction on $m$.  Note that $xp(x,y)=\sum_{i\geq 1,j \geq 0}a_{i-1,j}x^iy^j$.  Then
\begin{eqnarray*}
I_{k}(xp(x,y))&=&\sum_{i=0}^{k} \sum_{j=0}^{k-i} (-1)^j {k-i \choose j}a_{i-1,j} \\
&=&\sum_{i=1}^{k} \sum_{j=0}^{k-i} (-1)^j {k-i \choose j}a_{i-1,j} \mbox{ [since } a_{-1,j}=0 \mbox{ for all } j \mbox{]}\\
&=&\sum_{i'=0}^{k-1} \sum_{j=0}^{k-1-i'} (-1)^j {k-1-i' \choose j}a_{i',j}  \mbox{ [setting } i'=i-1\mbox{]}\\
&=&I_{k-1}(p(x,y)).
\end{eqnarray*}

For part (2), we first reverse the order of the sum to rewrite $I_k(p(x,y))=\sum_{j=0}^k \sum_{i=0}^{k-j}(-1)^j {k-i \choose j}a_{i,j}$, and note that
$yp(x,y)=\sum_{i\geq 0,j \geq 1}a_{i,j-1}x^iy^j$.  Then
\begin{eqnarray*}
I_{k}(yp(x,y))&=&\sum_{j=0}^{k} \sum_{i=0}^{k-j} (-1)^j {k-i \choose j}a_{i,j-1} \\
&=&\sum_{j=1}^{k} \sum_{i=0}^{k-j} (-1)^j {k-i \choose j}a_{i,j-1} \mbox{ [since } a_{i,-1}=0 \mbox{ for all } i \mbox{]}\\
&=&\sum_{j'=0}^{k-1} \sum_{i=0}^{k-j'-1} (-1)^{j'+1} {k-i \choose j'+1}a_{i,j'}  \mbox{ [setting } j'=j-1\mbox{]}\\
&=&(-1) \sum_{j'=0}^{k-1} \sum_{i=0}^{k-j'-1} (-1)^{j'} {k-i \choose j'+1}a_{i,j'}
\end{eqnarray*}

Now use the identity $\displaystyle{{k-i \choose j'+1}=\sum_{s=i}^{k-1}{s-i \choose j'} }$ to get 
$$I_k(yp)=(-1) \sum_{j'=0}^{k-1} \sum_{i=0}^{k-j'-1} (-1)^{j'}a_{i,j'} \sum_{s=i}^{k-1}{s-i \choose j'}. $$
For a fixed $s$ between $i$ and $k-1$ in this sum, we get 
$$(-1) \sum_{j'=0}^{k-1} \sum_{i=0}^{k-j'-1} (-1)^{j'}{s-i \choose j'}a_{i,j'} =(-1) \sum_{j'=0}^{s} \sum_{i=0}^{s-j'} (-1)^{j'} {s-i \choose j'}a_{i,j'} =(-1)I_s(p(x,y))$$
since terms where $j'>s$ and $i>s-j'$ give ${s-i \choose j'}=0$.  The result now follows by summing over $s$.

\qed

The proof of the next lemma uses Lemma~\ref{L:shift}. 

\begin{lem} \label{L:power} For $n, k \geq 0$, we have  $I_n(((x-1)(y-1))^k)=1$.
\end{lem}

\proof  We use induction on $n+k$.  Let $C_k=((x-1)(y-1))^k$.  The base cases of $n+k\leq 1$ are easy to check.  Then
\begin{eqnarray*}
I_n(C_k)&=&I_n((xy-x-y+1)C_{k-1}) \\
&=&I_n(xyC_{k-1})-I_n(xC_{k-1})-I_n(yC_{k-1})+I_n(C_{k-1}) \\
&=&I_{n-1}(yC_{k-1})-I_{n-1}(C_{k-1})-I_{n}(yC_{k-1})+I_n(C_{k-1}) \\
\end{eqnarray*}
by Lemma~\ref{L:shift}(1).

Now $I_n(C_{k-1})=I_{n-1}(C_{k-1})=1$ by induction.  For the remaining terms, we use Lemma~\ref{L:shift}(2).  
$$I_{n-1}(yC_{k-1})=(-1)\sum_{s=0}^{n-2} I_{s}(C_{k-1}) \mbox{ and } I_{n}(yC_{k-1})=(-1)\sum_{s=0}^{n-1} I_{s}(C_{k-1}).$$
Thus, by induction, we have $I_{n-1}(yC_{k-1})=1-n$ and $I_{n}(yC_{k-1})=-n$.  Hence $I_n(C_{k})=(1-n)-1-(-n)+1=1$, as desired.

\qed

Note that $\displaystyle{ ((x-1)(y-1))^{k}=\sum_{i=0}^k\sum_{j=0}^k (-1)^{i+j}{k\choose i}{k \choose j}x^iy^j}$.  Then Lemma~\ref{L:power} has the following purely combinatorial identity as a corollary.

\begin{cor}\label{C:iden}
For all $k, n \geq 0$,
$$\sum_{i=0}^n \sum_{j=0}^{n-i} (-1)^i {n-i \choose j} {k\choose i}{k \choose j}=1.$$
\end{cor}

It would be of interest to find a combinatorial proof of this identity.

 $I_k(p(x,y))$ can be computed as the trace of a matrix product.  Let $M_k$ be the $(r+1)\times (n+1)$ matrix with $(i,j)$ entry $(-1)^{j+1}{k-i+1 \choose j-1}$ and let $B$ be the $(n+1)\times (r+1)$ matrix of coefficients of the  polynomial $p(x,y)=\sum_{i,j \geq 0}b_{i,j}x^iy^j$, so the $(i,j)$ entry of the matrix $B$ is $b_{i-1,j-1}$.
$$M_k=
\left(
\renewcommand{\arraystretch}{1.5}
{\begin{array}{ccccc}
 {k \choose 0}& -{k \choose 1} & {k \choose 2} & \dots  \\
 {k-1 \choose 0}& -{k-1 \choose 1} & {k-1 \choose 2} & \dots  \\
  {k-2 \choose 0}&-{k-2 \choose 1} & {k-2 \choose 2} &  \dots  \\
\vdots & \vdots & \vdots & \vdots 
\end{array}}
\right)
\hskip.25in
B=
\left(
\renewcommand{\arraystretch}{1.5}
{\begin{array}{ccccc}
 b_{0,0} & b_{1,0} & b_{2,0} &  \dots  \\
 b_{0,1} & b_{1,1} & b_{2,1} & \dots   \\
 b_{0,2} & b_{1,2} & b_{2,2} &  \dots   \\
\vdots & \vdots & \vdots & \vdots 
\end{array}}
\right)$$

We omit the computational  proof of the next proposition.
\begin{prop}\label{P:trace}  Let $M_k$ be the $(r+1)\times (n+1)$ matrix with $(i,j)$ entry $(-1)^{i+j}{k-i+1 \choose j-1}$ and let $B$ be the $(n+1)\times (r+1)$ matrix of coefficients of the Tutte polynomial with $(i,j)$ entry $b_{i-1,j-1}$.  Then $I_k(T(G;x,y))=tr(M_k B)$.

\end{prop}

We write $r$ for $r(S)$ throughout the remainder of this section, and refer to $G$ as a rank $r$ ranked set.  Lemmas~\ref{L:k=r}, \ref{L:k>r} and \ref{L:r=n} all refer to the Tutte polynomial $T(G;x,y)=\sum b_{i,j}x^iy^j$ of a ranked set $G=(S,r)$. We first prove Theorem~\ref{T:affine} in the $k=r$ case. 

\begin{lem}\label{L:k=r}  ($r=k<n$ case)    Suppose $r=k$ and $r<n$.  Then
$$\sum_{i=0}^r \sum_{j=0}^{r-i} (-1)^j {r-i \choose j}b_{i,j}=0.$$
\end{lem}

\proof  (After Brylawski \cite{bry3}.)
Let $f(G;u,v)=\sum_{A \subseteq E}u^{r-r(A)}v^{|A|}$ be the {\it corank-cardinality} polynomial for $G$.  Then note that $$f(G;u,v)=v^{r}S\left(G;\frac{u}{v},v\right).$$

Using the relation $S(G;u,v)=T(G;u+1,v+1)$, we get $$f(G;u,v)=v^{r}T\left(G;\frac{u}{v}+1,v+1\right).$$  Setting $u=1$ gives 
$$(v+1)^n=f(G;1,v)=v^{r}T\left(G;\frac{v+1}{v},v+1\right).$$ 

Expanding the Tutte polynomial gives 
\begin{eqnarray*}
(v+1)^n&=&v^{r}T\left(G;\frac{v+1}{v},v+1\right) \\
&=&v^{r} \sum_{i, j \geq 0} b_{i,j} \left(\frac{v+1}{v}\right)^i (v+1)^j \\
&=& \sum_{i, j \geq 0} b_{i,j} v^{r-i}  (v+1)^{i+j}.
\end{eqnarray*}

Setting $y=v+1$ and expanding gives $$y^n=\sum_{i, j \geq 0} b_{i,j} (y-1)^{r-i} y^{i+j}= \sum_{i,j \geq 0} b_{i,j}\sum_{k=0}^{r-i} (-1)^k {r-i\choose k}y^{r+j-k}.$$

Thus, 
\begin{equation}\label{Eq:k=r}
y^{n-r}=\sum_{i,j \geq 0} b_{i,j}\sum_{k=0}^{r-i} (-1)^k {r-i\choose k}y^{j-k}.
\end{equation}

This is a Laurent polynomial  identity, so the constant term on the right-hand side (obtained by setting $k=j$) is identically zero since $r<n$.  This gives $$\sum_{i=0}^r \sum_{j=0}^{r-i} (-1)^j {r-i \choose j}b_{i,j}=0.$$
(The upper limit for the index $i$ can be taken as $r$ since $b_{i,j}=0$ for $i>r$.)

\qed

We  now use induction to prove the $k>r$ case of Theorem~\ref{T:affine}.

\begin{lem}\label{L:k>r} ($r<k<n$ case)    Suppose $r<k<n$.  Then
$$\sum_{i=0}^k \sum_{j=0}^{k-i} (-1)^j {k-i \choose j}b_{i,j}=0.$$
\end{lem}

\proof
We use induction on $k-r$.  If $k=r$, then the result is given by Lemma~\ref{L:k=r}.  Now assume $k>r$, set $m=k-r$ and the result is true for all $G$ with $k-r<m$.  Let $G'=G\times p$ be the free co-extension of $G$.  Then, from Prop.~\ref{P:exten}, we have  $G'/p=G$, so 
$$T(G;x,y)=T(G';x,y)-T(G'-p;x,y).$$
Since the $\rho(G')=\rho(G'-p)=\rho(G)+1$, we have $k-r=m-1$ for $G'-p$ and $G'$, so the result holds for both $G'$ and $G'-p$ by induction. The result now follows for $G$.
\qed

We will need one more lemma before we can complete the proof of Theorem~\ref{T:affine}. In our proof of Lemma~\ref{L:k=r}, we examined the constant term of  equation~(\ref{Eq:k=r}): $$y^{n-r}=\sum_{i,j \geq 0} b_{i,j}\sum_{k=0}^{r-i} (-1)^k {r-i\choose k}y^{j-k}.$$  We will examine non-constant terms from this equation to prove the next lemma.

\begin{lem}\label{L:r=n}  ($k<r=n$ case) Suppose $r=n$.  Then, for $k\leq n-1$, we have
$$\sum_{i=0}^k \sum_{j=0}^{k-i} (-1)^j {k-i \choose j}b_{i,j}=0.$$

\end{lem}
\proof
We use  equation~(\ref{Eq:k=r}) from Lemma~\ref{L:k=r}, setting $n=r$. It will be convenient to multiply both sides of this equation by $y^n$:

\begin{equation}\label{Eq:n=r}
y^n=\sum_{i,j \geq 0} b_{i,j}\sum_{k=0}^{n-i} (-1)^k {n-i\choose k}y^{n+j-k}.
\end{equation}

As in the proof of Lemma~\ref{L:k=r}, this is a polynomial identity. We examine the coefficients of $y^m$  for $0\leq m \leq n-1$.  Write 
$$\sum_{i,j \geq 0} b_{i,j}\sum_{k=0}^{n-i} (-1)^k {n-i\choose k}y^{n+j-k}=\sum_{m\geq 0} A_m y^m.$$

List the coefficients of $y^m$ for $0\leq m\leq n-1$ as a vector: $v_A:=\langle A_{0}, A_{1}, \dots, A_{n-1}\rangle$ and let $v_I:=\langle I_0, I_1, \dots, I_{n-1}\rangle$ be the vector of the $I_k$, where $$I_k=\sum_{i=0}^k \sum_{j=0}^{k-i} (-1)^j {k-i \choose j}b_{i,j}.$$

Now let $N$ be the $n \times n$ lower triangular  matrix with $N_{i,j}=(-1)^{n-i+1}{n-j \choose i-j}$.  Note that $N^{-1}=N$.  Then $N$  maps  $v_I^T$ to $v_A^T$, i.e., $Nv_I^T=v_A^T$. (This can be verified inductively, or by a direct expansion of equation~(\ref{Eq:n=r}).)

But, from equation~(\ref{Eq:n=r}), we have $A_m=0$ for all $0 \leq m \leq n-1$.  Thus $Nv_I=\overline{0}$, so $v_I=\overline{0}$, i.e., $I_k=0$ for all $0 \leq k \leq n-1$, as desired. 
\qed

We note that if $r=n$ for a matroid $M$, then the matroid consists of $n$ isthmuses, so $T(M;x,y)=x^n$, so Lemma~\ref{L:r=n} is trivial for matroids.

\begin{exmp}  We demonstrate  the matrix equation used in the proof of Lemma~\ref{L:k=r}  for $n=r=5$. Then we have the following relations:
\begin{eqnarray*} \hline \hline
A_{0}&=& -b_{0,0} \\
A_{1}&=& 5 b_{0,0}-b_{0,1}+b_{1,0} \\
 A_{2}&=&-10 b_{0,0}+5 b_{0,1}-b_{0,2}-4 b_{1,0}+b_{1,1}-b_{2,0} \\
A_{3}&=& 10 b_{0,0}-10 b_{0,1}+5 b_{0,2}-b_{0,3}+6 b_{1,0}-4 b_{1,1}+b_{1,2}+3 b_{2,0}-b_{2,1}+b_{3,0} \\
A_{4}&=& -5 b_{0,0}+10 b_{0,1}-10 b_{0,2}+5 b_{0,3}-b_{0,4}-4 b_{1,0}+6 b_{1,1}-4 b_{1,2}+b_{1,3}-3 b_{2,0} \\
& & +3 b_{2,1}-b_{2,2}-2 b_{3,0}+b_{3,1}-b_{4,0} \\ \hline \hline
\end{eqnarray*}

Since $Nv_I=v_A$ and $Nv_A=v_I$,  we can express the  $I_k$  as a linear combination of the coefficients $A_i$ for $i\leq k$, and vice versa.
$$
\left(
\begin{array}{ccccc}
 -1 & 0 & 0 & 0 & 0 \\
 4 & 1 & 0 & 0 & 0 \\
 -6 & -3 & -1 & 0 & 0 \\
 4 & 3 & 2 & 1 & 0 \\
 -1 & -1 & -1 & -1 & -1 \\
\end{array}
\right)
\left(
\begin{array}{c}
I_0 \\ I_1 \\ I_2 \\I_3 \\I_4
\end{array}
\right)=
\left(
\begin{array}{c}
A_0 \\ A_1 \\ A_2 \\A_3 \\A_4
\end{array}
\right)
$$

\end{exmp}

We now prove Theorem~\ref{T:affine}.

\proof [Proof Theorem~\ref{T:affine}.]
\begin{enumerate}
\item Assume $k<n$. We use a double induction, first on $k$, and, for fixed $k$, on $n-k$.  The result is trivial for $k=0$ and $n>0$.

Now let $k>0$ assume $I_m=0$ for all $m<k$ and, for fixed $m$, for  all $n>m$.  If $n-k=1$ (the base case for the induction on $n-k$), then $I_k=0$ by Lemma~\ref{L:k=r} (if $k=r$), Lemma~\ref{L:k>r} (if $k>r$) or Lemma~\ref{L:r=n} (if $k<r$, so $r=n$).  

Now assume $n-k>1$ and let $e \in S$.  By Theorem~\ref{T:tutte-gen}(1), we have
\begin{equation}\label{Eq:tutte}
T(G;x,y)=(x-1)^{r(G)-r(G-e)}T(G-e;x,y)+(y-1)^{1-r(e)}T(G/e;x,y).
\end{equation}

Expanding $(x-1)^{r(G)-r(G-e)}$ and $(y-1)^{1-r(e)}$ (if the exponents are positive), we can use Lemma~\ref{L:shift} to express $I_k$ as a linear combination of $I_m$ for various values of $m\leq k$ where $|S-e|= n-1$.  By induction on $k$ (when $m<k$) and $n-k$ (when $m=k$), the corresponding identities hold for each term on the right hand side of the recursion of equation~(\ref{Eq:tutte}), so they hold for $T(G;x,y)$.

\item When $n=k$, we use induction on $n-r$.
 Let $G'=(S,r')$ be formed from $G$ by changing the rank of $S$, leaving the ranks of all other subsets alone, i.e., $r'(A)=r(A)$ for $A \neq S$, and $r'(S)=n$.  Then it is a routine exercise to show
$$S(G';u,v)=u^{n-r}S(G;u,v)-(uv)^{n-r}+1.$$
Substituting $u=x-1$ and $v=y-1$ gives 
$$T(G';x,y)=(x-1)^{n-r}T(G;x,y)-((x-1)(y-1))^{n-r}+1.$$
For the left-hand side of this equation, we have $r(G')=n$, so $I_n(T(G'))=1$ (since this is the coefficient on the left-hand side of equation~(\ref{Eq:k=r}) used in Lemma~\ref{L:k=r}).  For the right-hand side, note that, by Lemma~\ref{L:shift}(1), we have 
$$I_n((x-1)^{n-r}T(G))=\sum_{i=0}^{n-r}(-1)^{n-r-i} {n-r \choose i} I_{n-i}(T(G)).$$
When $i>0$, we have $I_{n-i}(T(G))=0$ by part 1 of this theorem, so 
$$I_n((x-1)^{n-r}T(G))=(-1)^{n-r}I_n(T(G)).$$  Thus, 
$I_n(T(G))=(-1)^{n-r}I_n((x-1)(y-1)^{n-r})=(-1)^{n-r}$ by Lemma~\ref{L:power}.
\end{enumerate}

 \qed


\section{Antimatroids}\label{S:anti} Antimatroids are a class of greedoids of full rank whose feasible sets are closed under taking unions. See \cite{bz} for a treatment of antimatroids as a class of greedoids, or \cite{ej}, where antimatroids are used as a model for generalized convexity.

\begin{defn}\label{D:anti}  Let $G=(S,r)$ be a greedoid with rank function $r$ and feasible sets $\mathcal{F}$.  Then $G$ is an {\it antimatroid} if $r(S)=|S|$ and $F_1 \cup F_2 \in \mathcal{F}$ whenever $F_1, F_2 \in \mathcal{F}$.
\end{defn}
A set is {\it convex} if its complement is feasible. Then we have an expansion of the Tutte polynomial in terms of the convex sets. For $C$ convex, let $ex(C)$ be the {\it extreme} points of $C$, i.e., $p \in ex(C)$ if $p \in C$ but $p \notin \overline{C-p}$, where $\overline{A}$ is the {\it convex closure} of $A$, i.e., the smallest convex set containing $A$. We write $int(C):=C-ex(C)$ for the interior of $C$. The interior $int(C)$ has combinatorial interpretations for all the classes of antimatroids we consider here.

\begin{thm}[Theorem 2.2 \cite{gposet}] \label{T:anti}
Let $G$ be an antimatroid with convex sets $\mathcal{C}$. Then
$$T(G;x,y)=\sum_{C \in \mathcal{C}} (x-1)^{|C|}y^{|int(C)|}.$$
\end{thm}

 Let $a_{i,j}$ be the number of convex sets with $i$ points, $j$ of which are interior. Expanding $(x-1)^{|C|}$ allows us to express the Tutte coefficients $b_{i,j}$ in terms of the $a_{i,j}$.
 
 \begin{lem}\label{L:anti} Let $G=(S,r)$ be an antimatroid with $|S|=n$ and Tutte polynomial $T(G;x,y)=\sum b_{i,j}x^iy^j.$ Then $$b_{i,j}=\sum_{s=i}^n(-1)^{s-i}{s \choose i}a_{s,j}.$$
  \end{lem}
  
  Combining Lemma~\ref{L:anti} with Theorem~\ref{T:affine} gives us the following identities that all antimatroids satisfy.
  
  \begin{cor}\label{C:anti}
  Let $G=(S,r)$ be an antimatroid with $a_{i,j}$ convex sets of size $i$ and interior of size $j$.
\begin{enumerate}
\item  For $k<n$,
 $$\sum_{i=0}^k \sum_{j=0}^{k-i} (-1)^j {k-i \choose j}\sum_{s=i}^n(-1)^{s-i}{s \choose i}a_{s,j}=0.$$
 \item For $k=n$,
  $$\sum_{i=0}^n \sum_{j=0}^{n-i} (-1)^j {n-i \choose j}\sum_{s=i}^n(-1)^{s-i}{s \choose i}a_{s,j}=1.$$
\end{enumerate}
  \end{cor}
  
 The $k=0, 1$ and $2$ cases of Cor.~\ref{C:anti} (1) are worth separate consideration. 
 \begin{itemize}
\item {\bf $k=0.$} We say a convex set is {\it free} if it has empty interior. Let $f_i$ be the number of free convex sets with $i$ points, so $f_i=a_{i,0}$ (and assume $f_0=1$). Then for $k=0$, the identity from Cor.~\ref{C:anti} (1) reduces to 
  $$\sum_{i=0}^n (-1)^if_i =0.$$
This identity is Theorem 4.5 in \cite{ej}; Edelman and Jamison attribute this result to Lawrence. An equivalent formulation using the characteristic polynomial appears as Prop. 7 of \cite{gm3}.
\item {\bf $k=1.$} This identity gives $b_{1,0}=b_{0,1}$, i.e., the coefficient of $x$ equals the coefficient of $y$ in the Tutte polynomial. Translating to convex sets, we get:
\begin{equation}\label{eq:k=1}
\sum_{i=0}^n(-1)^{i-1}if_i=\sum_{i=0}^n(-1)^ia_{i,1}
\end{equation}

The invariant on the left-hand side is $b_{1,0}$, the beta invariant, and the right-hand side is an alternating sum over convex sets with exactly one interior point. The beta invariant gives interesting combinatorial information about the antimatroid -- this is the focus of \cite{gbeta}. We will examine this identity for several classes of antimatroids below.

We can rewrite the right hand side as a (double) sum over all interior points in $G$. 

\begin{cor}\label{C:b01}For $p\in G$, let $\mathcal{C}_p$ be the collection of convex sets in $G$ with unique interior point $p$. Then 
$$b_{0,1}=\sum_{i=0}^n(-1)^ia_{i,1}=\sum_{p \in int(G)} \sum_{C \in \mathcal{C}_p}(-1)^{|C|}$$
\end{cor}

This allows us to express the beta invariant as a sum over interior points of $G$. 
\item {\bf $k=2$} We use $b_{0,2}+b_{2,0}=b_{1,0}+b_{1,1}$. Simplifying  gives
$$\sum_{i=0}^n (-1)^i \left({i+1 \choose 2}f_i +ia_{i,1}+ a_{i,2}\right) =0.$$
There are several equivalent formulations. This identity, and those involving higher indices, give more involved combinatorial results for the antimatroids we treat below.
\end{itemize}

Our immediate goal is to interpret the $k=1$ identity for four families of antimatroids:  trees, posets, chordal graphs and finite subsets of $\mathbb{R}^n$. For each class, we describe the antimatroid structure by specifying the convex sets.

\subsection{Trees} Let $T$ be a  tree with edges $E$. The {\it pruning antimatroid} $G=(E,r)$ is defined on the set of edges of $T$, where the   convex sets are the subtrees of $T$. An edge is {\it interior} if it is not a leaf of the subtree.

Then a convex set (subtree) is {\it free} if it has no interior edges, i.e., the edges form a star. We now give a combinatorial interpretation to the beta invariant for a tree. We omit the proof, which follows from Prop. 1.7 and Cor. 3.3 of \cite{ganti} (or can be proven directly).

\begin{prop}\label{P:treebeta}
Let $T$ be a tree with $n$ edges and $m$ interior edges. Then $$\sum_{i=0}^n (-1)^{i}if_i=m.$$
\end{prop}

Then the $k=1$ identity translates to the following theorem.  

\begin{thm}\label{T:treek=1} 
Let $T$ be a tree with $m$ interior edges and let $\mathcal{S}$ be the collection of all subtrees with exactly one interior edge. Then
$$\sum_{S \in \mathcal{S}}(-1)^{|S|}=-m.$$
\end{thm}

\begin{figure}[h]
\begin{center}
\includegraphics[width=3in]{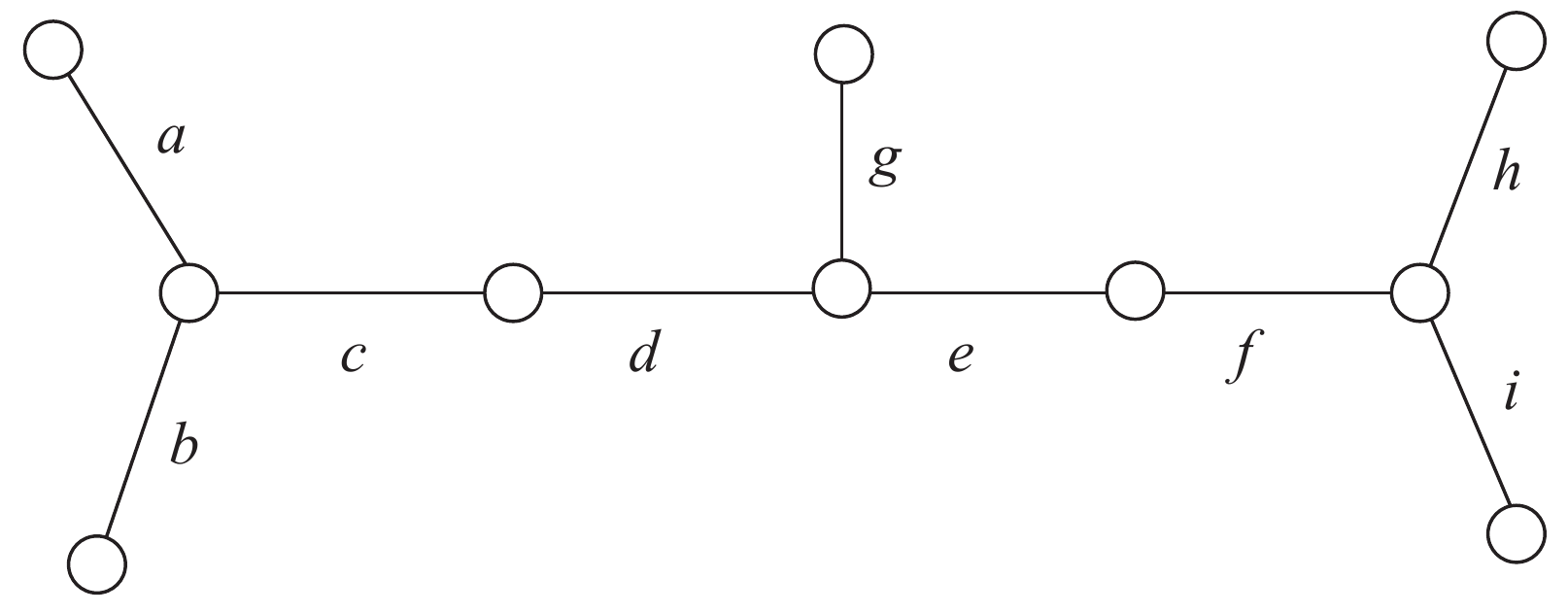}

\caption{A tree.}\label{F:tree}
\end{center}
\end{figure}

As an example, consider the tree in Fig.~\ref{F:tree}. There are nine subtrees of size 3 with exactly one interior edge, six of size 4, and one of size 5. Then Theorem~\ref{T:treek=1} gives $\sum_{S \in \mathcal{S}}(-1)^{|S|}=-9+6-1=-4,$ as required.

When $k=2$, the identity involves $f_i (=a_{i,0}), a_{i,1}$ and $a_{i,2}$. We list the values for these invariants for the tree of Fig.~\ref{F:tree} in Table~\ref{Ta:tree data}.
\begin{table}[htdp]
\begin{center}
\begin{tabular}{c || c c c c c c c}
 & $i=0$ & $i=1$ & $i=2$ & $i=3$ & $i=4$ & $i=5$ & $i=6$ \\ \hline
 $f_i=a_{i,0}$ & 1 & 9 & 11 & 3 & 0 & 0 & 0\\
  $a_{i,1}$ & 0 & 0 & 0 & 9 & 6 & 1 & 0 \\
   $a_{i,2}$ & 0 & 0 & 0 & 0 & 6 & 5 & 1 \\
\end{tabular}
\bigskip
\caption{Data for the tree of Fig.~\ref{F:tree}.}
\label{Ta:tree data}
\end{center}
\end{table}%

Suppressing 0-terms, the $k=2$ identity becomes 
$$(-f_1+3f_2-6f_3)+(-3a_{3,1}+4a_{4,1}-5a_{5,1})+(a_{4,2}-a_{5,2}+a_{6,2})$$
This reduces to $(-9+33-18)+(-27+24-5)+(6-5+1)=0$, as required.
\subsection{Posets} There are several ways to give a poset $P$ an antimatroid structure. We use the {\it double shelling antimatroid}. Recall a set $I$ is an {\it order ideal} in $P$ if $x \in I$ and $y\leq x$ implies $y \in I$. $J$ is an {\it order filter} if
$x \in J$ and $y\geq x$ implies $y \in J$. Then a set $C$ is convex in the double shelling antimatroid associated with $P$  if there is an order ideal $I$ and an order filter $J$ such that $x \in C$ precisely when $x\geq y$ for all $y\in I$ and $x \leq z$ for all $z \in J$.

A convex set $C$ is {\it free} if it contains no chains of length greater than 2, i.e., every element of $C$ is either minimal or maximal in the poset $P-I-J$. A {\it bottleneck} in a poset is an element that is not maximal or minimal, but is comparable to every element of the poset. Then the following  combinatorial characterization of the beta invariant was proven by Edelman and Reiner in \cite{er}.

\begin{thm}[Cor. 4.4 \cite{er}]\label{T:betaposet}
Let $P$ be a poset with $f_i$ free convex sets and $b$ bottlenecks. Then 
$$\sum_{i=0}^n (-1)^{i}if_i=b.$$
\end{thm}

To interpret the $k=1$ identity here, we note that a convex set $C$ has exactly one interior point if $C$ is convex and there is a unique $x\in C$ with $x \notin ex(C)$. Then there is a chain of length 3 in $C$, with $x$ in the middle.

\begin{thm}\label{T:posetk=1} Let $P$ be a poset with $b$ bottle necks and let $\mathcal{S}$ denote the set of all convex sets with exactly one interior point. Then
$$\sum_{S \in \mathcal{S}} (-1)^{|S|}=-b.$$
\end{thm}

\begin{figure}[h]
\begin{center}
\includegraphics[width=1in]{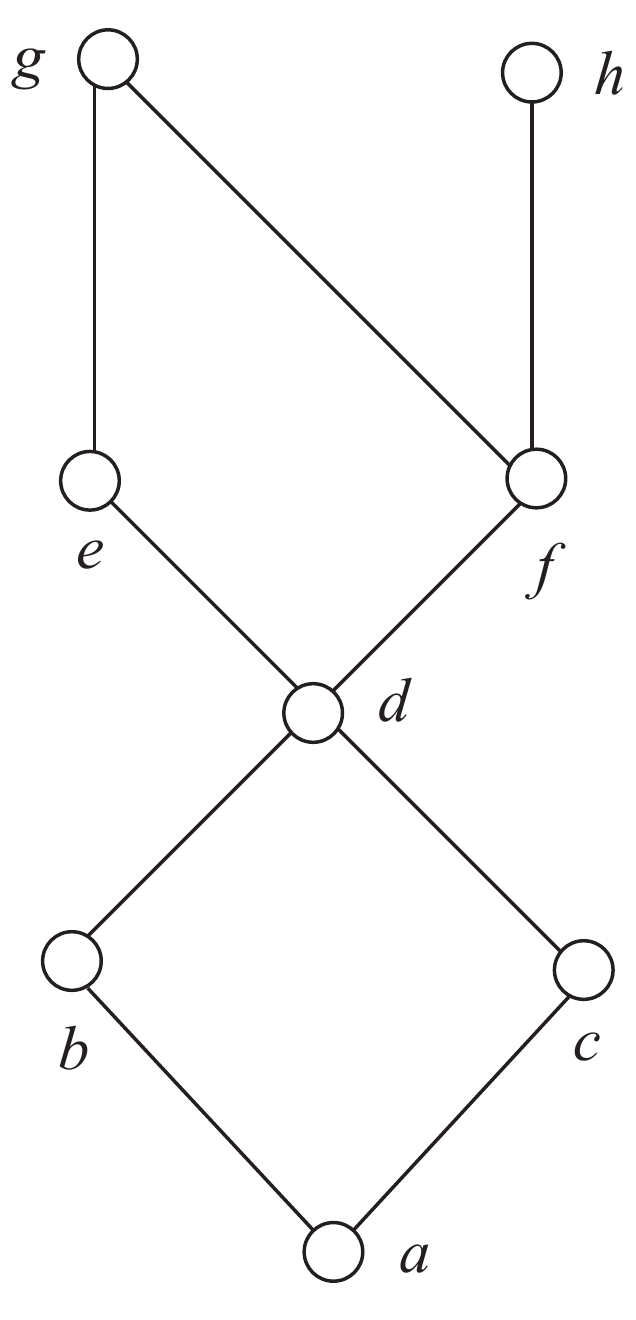}

\caption{A poset.}\label{F:poset}
\end{center}
\end{figure}

As an example, consider the poset of Fig.~\ref{F:poset}. The element $d$ is a bottleneck. We list the convex sets with exactly one interior point in Table~\ref{Ta:poset}.
\begin{table}[h]
\begin{center}
\begin{tabular}{l | l l l l l }
Size & & & & & \\ \hline
3 & $bde$ & $bdf$ & $cde$ & $cdf$ & $dfh$ \\
4 & $bcde$ & $bcdf$ & $bdef$ & $cdef$ & $defh$ \\
5 & $bcdef$ & & & 
\end{tabular}
\bigskip
\caption{The convex sets with exactly one interior element.}
\label{Ta:poset}
\end{center}
\end{table}%
The alternating sum then gives $-4+4-1=-1$, as required. 

We leave full consideration of the $k=2$ case to the interested reader. Some care must be taken in determining the convex sets with exactly two interior elements, however. For instance, in the poset of Fig.~\ref{F:poset}, the chain $b<d<f<h$ is a convex set with two interior elements, but the chain $a<b<d<e$ is not, since $abde$ is not a convex set in this antimatroid.

\subsection{Chordal graphs} When $G$ is a chordal graph, a vertex is {\it simplicial} if its neighbors form a clique. The ground set of the {\it simplicial shelling antimatroid} is the collection of vertices of $G$. An ordered set $\{v_1, v_2, \dots, v_k\}$ is {\it shellable} if $v_1$ is simplicial, $v_2$ is simplicial in $G-v_1$, $v_3$ is simplicial in $G-\{v_1,v_2\}$, and so on. Then we can remove the dependence on an order by defining a set to be feasible if there is some shellable ordering of the vertices. As usual, the convex sets are the complements of the feasible sets.

A set of vertices is free convex if the graph it induces is a clique in $G$. The next result gives a combinatorial interpretation to the beta invariant $b_{1,0}.$
\begin{thm} [Theorem 5.1 \cite{gbeta}] \label{T:chordalbeta} Let $G$ be a chordal graph with $b$ 2-connected blocks, and let $f_i$ be the number of cliques of size $i$. Then $$\sum_{i=0}^n (-1)^{i}if_i=b-1.$$
\end{thm}

\begin{figure}[h]
\begin{center}
\includegraphics[width=2.5in]{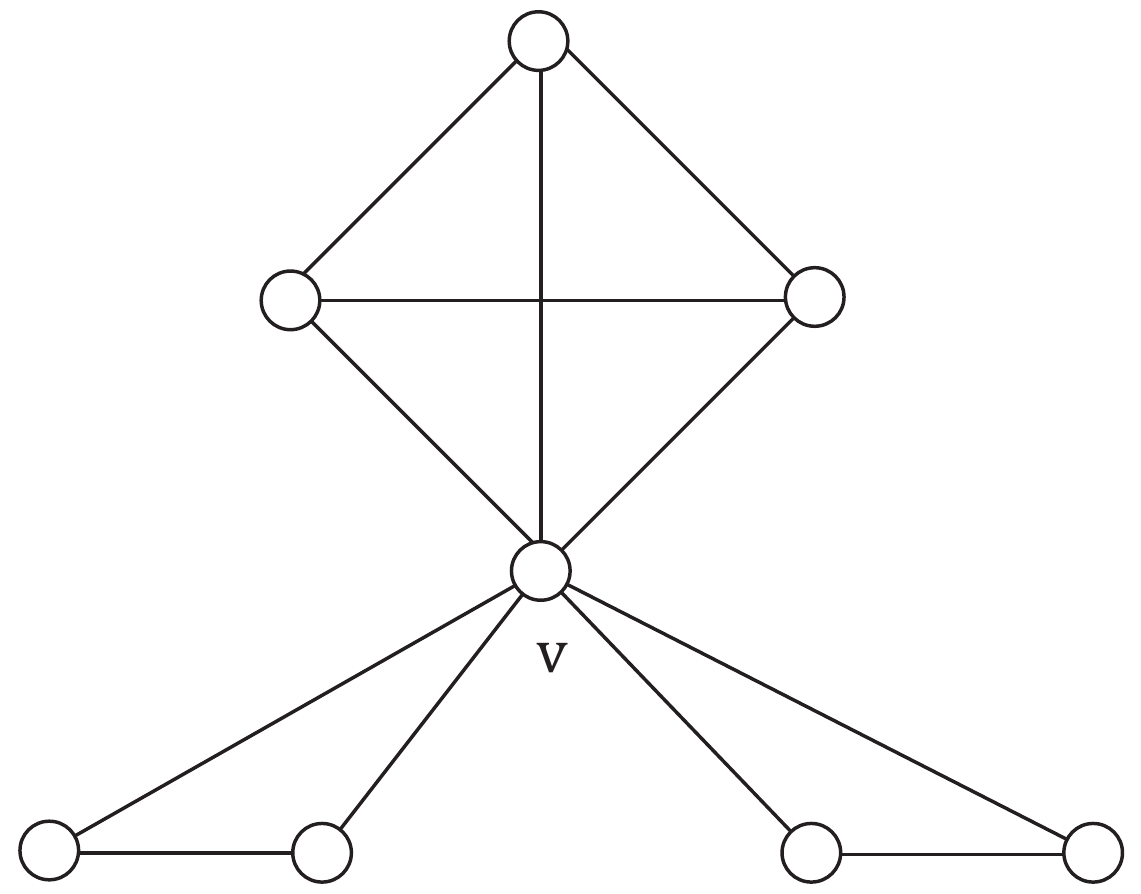}

\caption{A 1-sum of cliques: The vertex $v$ is the unique internal vertex in the convex set in the chordal graph.}
\label{F:blocksum}
\end{center}
\end{figure}

To interpret the $k=1$ identity for chordal graphs, we need the following characterization of convex sets with a single interior vertex. A {\it 1-sum} of cliques is formed by gluing a disjoint collection of cliques together at  one common vertex, as in Fig.~\ref{F:blocksum}. These are the convex sets we need.

\begin{lem}\label{L:chordal}
Let $G$ be a chordal graph. A subset $C \subseteq V$ is convex with exactly one interior vertex if and only if the subgraph induced by $C$ is a 1-sum of cliques.
\end{lem}

Then the $k=1$ identity takes the following form for chordal graphs.

\begin{thm}\label{T:chordalk=1}
Let $G$ be a chordal graph and with $b$ blocks and let $\mathcal{S}$ denote the set of all convex sets with exactly one interior vertex. Then
$$\sum_{S \in \mathcal{S}} (-1)^{|S|}=1-b.$$
\end{thm}

\begin{figure}[h]
\begin{center}
\includegraphics[width=2.75in]{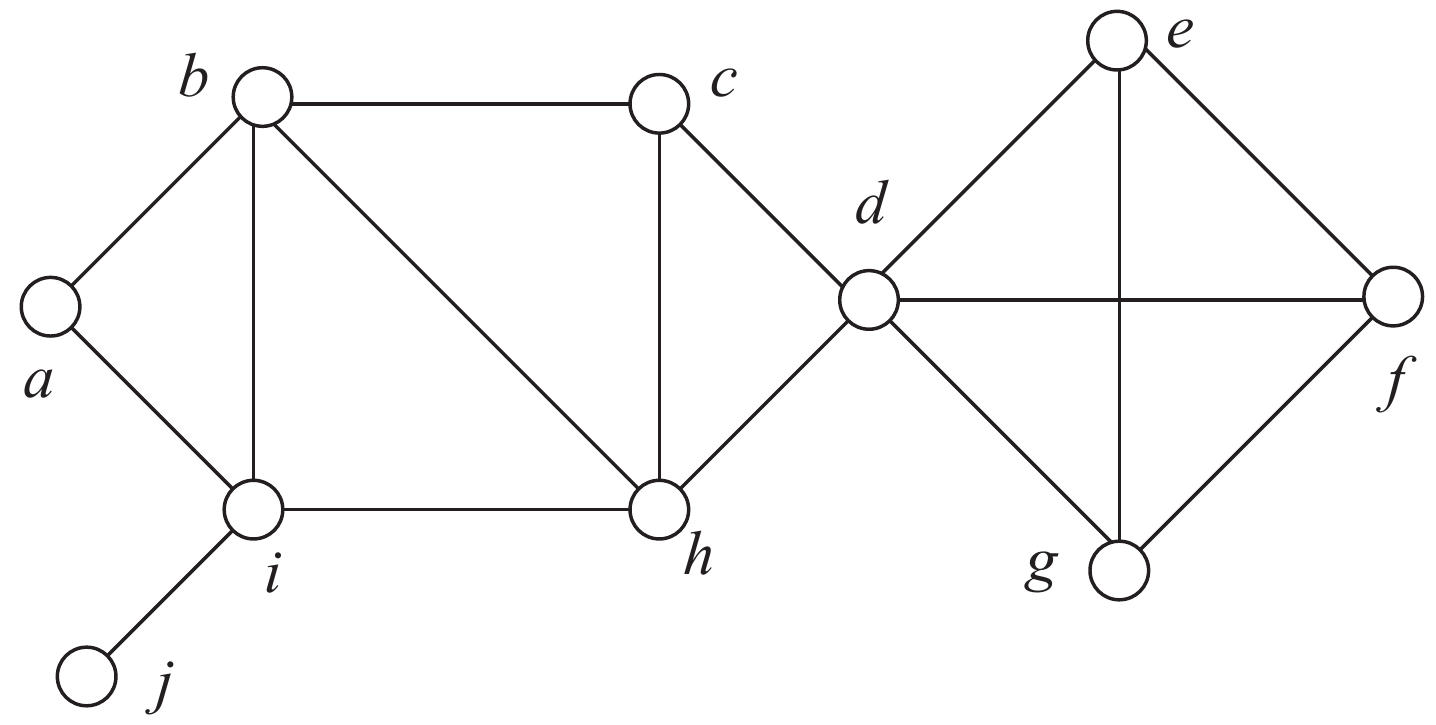}

\caption{A chordal graph.}
\label{F:chordal}
\end{center}
\end{figure}

As an example, consider the chordal graph of Fig.~\ref{F:chordal}. Note that $G$ has 3 blocks. Then  $d$ and $i$ are the only vertices that are unique interior vertices in any convex set. We compute $b_{0,1}$ using Cor.~\ref{C:b01}. Then we see the vertex $d$ is in six such sets of size 3, nine of size 4, five of size 5, and one of size 6. For the vertex $i$, the corresponding counts give three sets of size 3 and two of size 4. Then the alternating sum from Theorem~\ref{T:chordalk=1} gives us $-9+11-5+1=-2,$ as required. As with posets, some care is needed to ensure a subset of vertices gives a convex set. For instance, although the vertex set $\{d,h,i\}$ gives an induced subgraph that corresponds to a 1-sum of blocks, the set $\{d,h,i\}$ is not convex.

Trees are chordal, so Theorem~\ref{T:chordalk=1} applies to the pruning greedoid associated to the vertices of a tree.  In this context, the only free, convex sets are the single vertices and  pairs of adjacent vertices. A subset of vertices has precisely one interior vertex if it corresponds to a star, where the central vertex is the interior point. (Note that this agrees with our more general interpretation given in Lemma~\ref{L:chordal}.) 

Thus, a vertex $v$ of degree $d$ will contribute $\sum_{i=2}^d (-1)^{i-1} {d \choose i}$ to the coefficient $b_{0,1}$.  Since $\sum_{i=2}^d (-1)^{i-1} {d \choose i}=1-d$, it is easy to see that summing over all convex sets with a unique interior point gives $2-n$, which agrees with Theorem~\ref{T:chordalk=1} since there are $n-1$ blocks (the edges of $T$).

\subsection{Finite point sets in Euclidean space} Finite subsets of $\mathbb{R}^n$ are the prototypical examples of convex geometries, or, dually, antimatroids. A subset $C$ of a finite set $S$ is {\it convex} if $C=\overline{C} \cap S$, where $\overline{C}$ is the convex hull of $C$ in $\mathbb{R}^n$. In this context, interior points are straightforward. As usual, the complements of the convex sets are the feasible sets of the antimatroid.

The beta invariant counts the number of interior points, with a sign indicating the parity of the dimension of the set $S$. This theorem is the main result of \cite{er}.

\begin{thm}[Theorem 1.1 \cite{er}] \label{T:betapoints} 
Let $S$ be a finite subset of $\mathbb{R}^n$. Then 
$$\sum_{i=0}^n (-1)^i i f_i =(-1)^{n-1}|int(S)|.$$
\end{thm}

Combining Theorem~\ref{T:betapoints} with the $k=1$ identity gives us the next result.

\begin{cor}\label{C:pointset}
Let $S$ be a finite subset of $\mathbb{R}^n$. Let $\mathcal{C}_1$ be the collection of all convex sets with exactly one interior point. Then
$$\sum_{C\in\mathcal{C}_1} (-1)^{|C|} =(-1)^{n}|int(S)|.$$
\end{cor}

\begin{figure}[h]
\begin{center}
\includegraphics[width=1.75in]{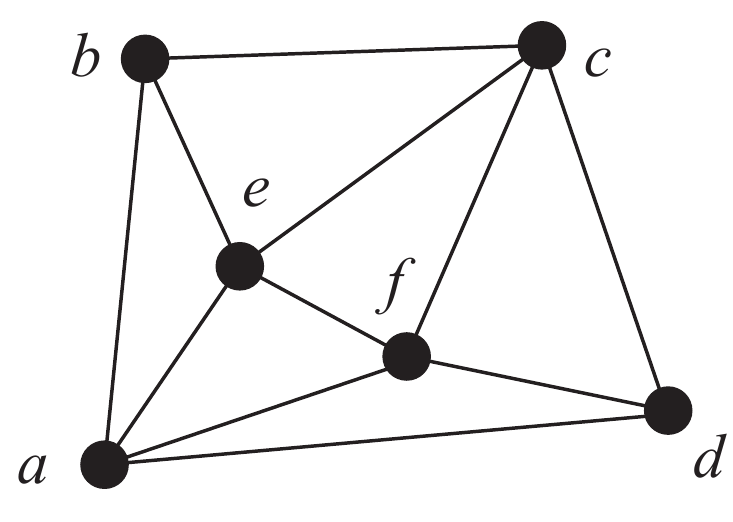}

\caption{A finite subset of the plane.}
\label{F:ptset}
\end{center}
\end{figure}

For example, the set $S$ of six points in the plane in Fig.~\ref{F:ptset} has 2 interior points, $e$ and $f$. (The point set is triangulated to help visualize the convex sets.)

Using the interpretation for $b_{0,1}$ given in Cor.~\ref{C:b01}, we find all convex sets that have $e$ as unique interior point, and those that have $f$. In this example, we have $abce, abef$ and $abcef$ are  the  convex sets with $e$ as the unique interior point, while $f$ is the unique interior point for $acdf, adef$ and $acdef$. This gives $b_{0,1}=4(-1)^4+2(-1)^5=2=|int(S)|$.

For finite point sets, it is easier to interpret the  identities for larger values of $k$. For instance,  the $k=2$ identity
depends on convex sets with 0, 1 or 2 interior points. In our example, we have the following data (see Table~\ref{Ta:point data}):

\begin{table}[htdp]
\begin{center}
\begin{tabular}{c || c c c c c c c}
 & $i=0$ & $i=1$ & $i=2$ & $i=3$ & $i=4$ & $i=5$ & $i=6$ \\ \hline
 $f_i=a_{i,0}$ & 1 & 6 & 15 & 15 & 6 & 1 & 0\\
  $a_{i,1}$ & 0 & 0 & 0 & 0 & 4 & 2 & 0 \\
   $a_{i,2}$ & 0 & 0 & 0 & 0 & 0 & 1 & 1 \\
\end{tabular}
\bigskip
\caption{Data for the point set of Fig.~\ref{F:ptset}.}
\label{Ta:point data}
\end{center}
\end{table}%
Then $\sum_{i=0}^n (-1)^i \left({i+1 \choose 2}f_i +ia_{i,1}+ a_{i,2}\right) =-(6)+(45)-(90)+(60+16)-(15+10+1)+(1)=0.$

\section{Examples and counterexamples}\label{S:ex}
In this section, we present a series of examples to show how the identities of Theorem~\ref{T:affine} can be used to help find the Tutte polynomial. Since matroids and greedoids (and, more generally, ranked sets) satisfy the same families of identities, it is natural to ask if the Tutte polynomial can distinguish these objects. Counterexamples are given here to show that a matroid and a greedoid (that is not a matroid) can share the same Tutte polynomial. Variations that use deletion and contraction are also given.

We begin by considering other relations the Tutte polynomial coefficients satisfy when the underlying object is a matroid. Here is the complete list of affine relations that form a basis for all affine relations satisfied by these coefficients.

\begin{thm} [Brylawski] \label{T:basis} Let $M$ be a matroid with no isthmuses and let $T(M;x,y)=\sum b_{i,j}x^iy^j$ be its Tutte polynomial.  Let $r(M)=r$ and $|S|=n$.  Then
\begin{enumerate}
\item $b_{i,j}=0$ for all $i>r$ and all $j>0$;
\item $b_{r,0}=1; b_{r,j}=0$ for all $j>0$;
\item $b_{r-1,0}=n-r; b_{r-1,j}=0$ for all $j>0$;
\item $b_{i,j}=0$ for all $i$ and $j$ such that $1 \leq i \leq r-2$ and $j \geq n-r$;
\item $b_{0,n-r}=1; b_{0,j}=0$ for all $j>n-r$;
\item $\displaystyle{\sum_{i=0}^k \sum_{j=0}^{k-i} (-1)^j {k-i \choose j}b_{i,j}=(-1)^{n-r}}$ for all $0 \leq k \leq n-1$.
\end{enumerate}
Further, these identities form a basis for all affine relations satisfied by the Tutte polynomial coefficients.
\end{thm}

Since these identities form a basis for all affine relations satisfied by the Tutte polynomial coefficients, they also generate the ``new'' relation of Corollary~\ref{C:matroid}.  We do not believe this identity has been explicitly stated before, however.

\begin{cor}\label{C:matroid} Let $M$ be a rank $r$ matroid on $n$ points.  Then if $T(M)=\sum b_{i,j}x^iy^j$, 
$$\sum_{i=0}^n \sum_{j=0}^{n-i} (-1)^j {n-i \choose j}b_{i,j}=(-1)^{n-r}.$$
\end{cor}

As an example of the $k=n$ identity, consider the two uniform matroids $U_{2,4}$ and $U_{3,4}$.
\begin{exmp}\label{E:uniform} Let $M$ be the uniform matroid $U_{2,4}$.  Then $T(M;x,y)=x^2+2x+2y+y^2$.  Recall 
\begin{eqnarray*}
 I_4(T(M))&=& b_{0,0}-4 b_{0,1}+6 b_{0,2}-4 b_{0,3}+b_{0,4}+b_{1,0}-3 b_{1,1}+3 b_{1,2}-b_{1,3} \\
 & & +b_{2,0}  -2 b_{2,1}+b_{2,2}+b_{3,0}-b_{3,1}+b_{4,0}.
\end{eqnarray*}
 Then  $I_4(T(M))=-8+6+2+1=1,$ as required by the theorem since $n-r$ is even.

For $M=U_{3,4}$, we have $T(M)=x^3+x^2+x+y$.  This time, we find $I_4(T(M))=-4+1+1+1=-1$, since $n-r$ is odd.

\end{exmp}

The next set of examples examines three different rank 3, cardinality 5 ranked sets (one matroid and two greedoids).

\begin{exmp}\label{Ex:basis}  Let $M$ be a rank 3 matroid on 5 points.  Using all 6 of the relations of Theorem~\ref{T:basis} forces $T(M)=x^3+2x^2+b_{1,0}x+b_{1,1}xy+b_{0,1}y+y^2$, where $b_{1,0}, b_{0,1}$ and $b_{1,1}$ are undetermined.  Then the simplified relations $I_k=0$ for $k \leq 4$ are:
\begin{itemize}
\item $I_0:  b_{0,0}=0$.
\item $I_1:  b_{0,1}=b_{1,0}$.
\item $I_2:  b_{1,1}+b_{1,0}=3$.
\item $I_3:  2b_{1,1}+2b_{1,0}=6$.
\item $I_4:  3b_{1,1}+3b_{1,0}=9$.
\end{itemize}
(Note that $I_2, I_3$ and $I_4$ all give equivalent relations.)  Now $I_5(T)=13-4b_{1,0}-4b_{1,1}$.  Using $I_2$, this simplifies to $I_5(T)=1$, as required by Cor.~\ref{C:matroid}.  Thus, the ``new'' relation given is determined by the relations of Theorem~\ref{T:basis}.

Now consider the three graphs $G, G'$ and $G''$ of Fig.~\ref{F:3graphs}.  Each has rank 3 on ground sets of size 5.  We give the Tutte polynomials of each, where $G$ is an unrooted graph, but $G'$ and $G''$ are rooted graphs, and so use the branching greedoid rank function.

\begin{figure}[h]
\begin{center}
$$G \hskip1.35in G' \hskip1.35in G''$$

\includegraphics[width=1in]{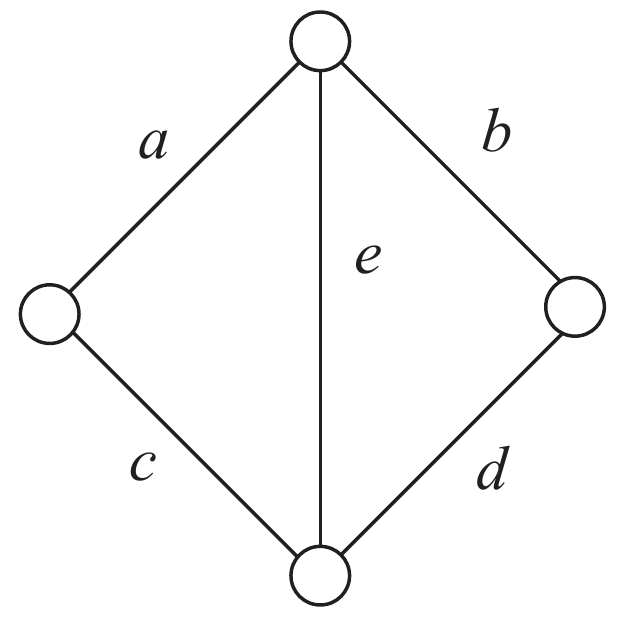} \hskip.4in 
\includegraphics[width=1in]{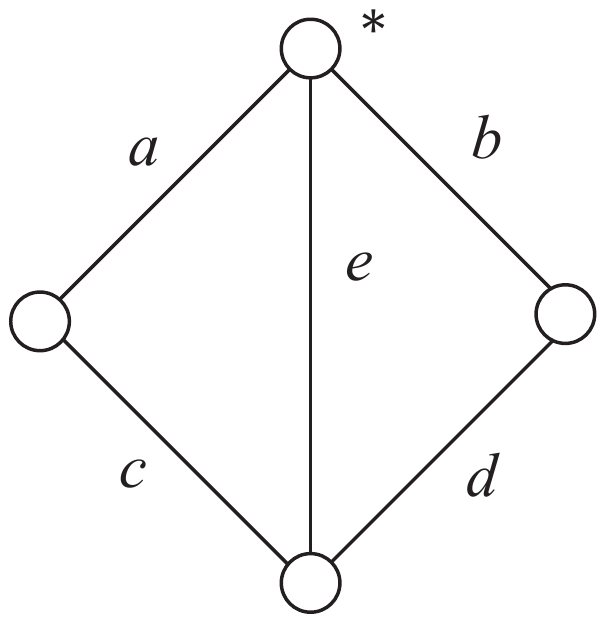} \hskip.4in
\includegraphics[width=1in]{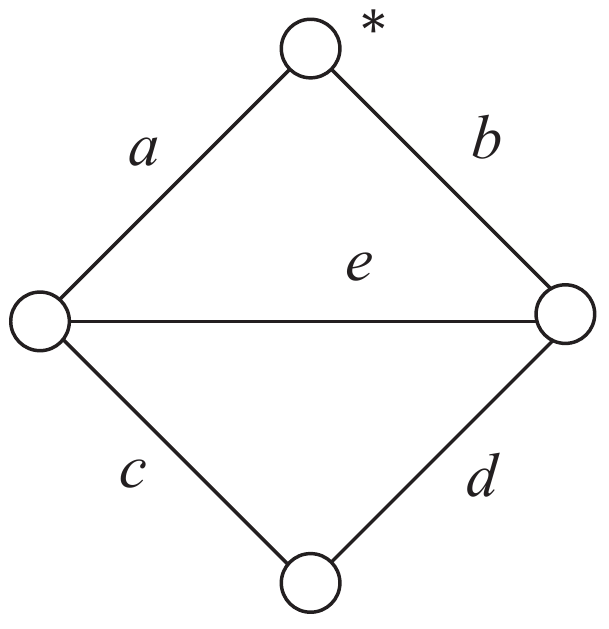}

\caption{One unrooted and two rooted graphs for Example~\ref{Ex:basis}.}
\label{F:3graphs}
\end{center}
\end{figure}

Then
\begin{eqnarray*}
T(G)&=&x^3+2 x^2+2 x y+x+y^2+y \\
T(G')&=&x^3 y^2-3 x^2 y^2+2 x^2 y+x^2+3 x y^2-2 x y+3 x+3 y\\
T(G'')&=&x^3 y^3-3 x^2 y^3+2 x^2 y+3 x y^3-3 x y+4 x-y^3+y^2+4 y
\end{eqnarray*}

Then we list  the various identities $I_0, \dots, I_5$ for each of these polynomials.  Note that we have used $I_0$ to simplify $I_1$, $I_0$ and $I_1$ to simplify $I_2$, and so on.  

\begin{table}[htdp]
\begin{center}
\begin{tabular}{|c|c|c|c|c|} \hline
$k$& $I_k$&$ T(G)$ & $T(G')$ & $T(G'')$ \\ \hline
0&$b_{0,0}$ & 0 & 0 & 0 \\ \hline
1&$b_{0,1}-b_{1,0}$ & $1-1$ & $3-3$ & $4-4$ \\ \hline
2&$b_{2,0}+b_{0,2}-b_{1,0}-b_{1,1}$ & $2+1-2-1$ & $1+0-(-2)-3$ & $0+1-(-3)+4$\\ \hline
3&$b_{3,0}-b_{2,1}+b_{1,2}$ & $1-0+0$ & $0-2+3$ & $0-2+0$\\ 
 & $-b_{0,3}-b_{2,0}+b_{0,2}$ &$-0-2+1$ & $0-1+0 $& $-(-1)-0+1$ \\ \hline
\end{tabular}
\bigskip
\caption{Calculations of the identities $I_k$ for the three examples of Example~\ref{Ex:basis}.}
\label{Ta:iden}
\end{center}
\end{table}
Note that the simplified matroid identity $I_2:  b_{1,1}+b_{1,0}=3$ is false for $T(G')$ and $T(G'')$ (as $G'$ and $G''$ are not matroids).
\end{exmp}

Part (6) of Theorem~\ref{T:basis} is the true in our more general setting; this is Theorem~\ref{T:affine}(1).  Which of the other parts of  \ref{T:basis} remain true?  In their full generality, only part (1) is still valid.

\begin{prop}\label{P:rest}  In addition to (6), (1) still holds:  $b_{i,j}=0$ for all $i>r$ and all $j>0$.
\end{prop}

\begin{exmp}\label{E:badiden}
Let $G$ is an antimatroid with ground set $\{1,2,3\}$ with feasible sets $$\{\emptyset,\{1\},\{3\},\{1,2\},\{1,3\},\{2,3\},\{1,2,3\}\}.$$ Then $T(G;x,y)=x^3 y-3 x^2 y+2 x^2+3 x y-x-y.$  Then $r(G)=|S|=3$, and $G$ has no isthmuses.  We note that $G$ is the edge pruning greedoid associated with a path on three edges. Note that 
\begin{itemize}
\item $b_{3,0}=0$  and $b_{3,1}=1$, so Theorem~\ref{T:basis}(2) is no longer valid.
\item $b_{2,0}=2$ and $b_{2,1}=-3$,  so Theorem~\ref{T:basis}(3) is also false.
\item $b_{1,1}=3$, so Theorem~\ref{T:basis}(4) is false, too.
\end{itemize}
\end{exmp}

When $G$ is not a matroid, we can have Tutte polynomials with negative coefficients.  But this property does not distinguish the class of greedoids from matroids. The next three counterexamples examine this limitation.

\begin{exmp}\label{Ex:pos}  Let $G$ be a greedoid with feasible sets $$\mathcal{F}=\{\emptyset,\{a\},\{b\},\{c\},\{a,b\},\{a,c\},\{b,c\},\{a,d\}\}.$$  Note that $\{d\}$ is not  feasible, but $\{a,d\}$ is, so $G$ is not a matroid.  This is the second truncation of the branching greedoid associated to the rooted tree in Figure~\ref{F:rtree}.

\begin{figure}[h]
\centerline{\includegraphics[width=1.35in]{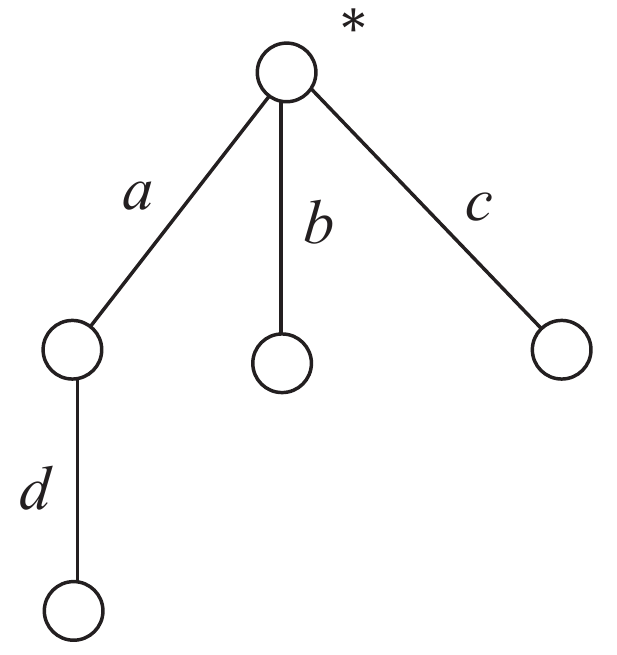}} \hskip.4in 

\centerline{\includegraphics[width=2.5in]{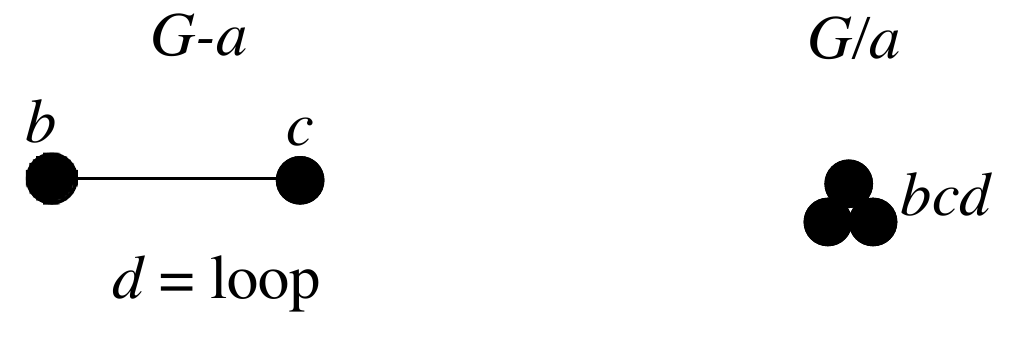}}
\caption{$G$ is the second truncation of the rooted tree above. $G-a$ and $G/a$ are both matroids.}
\label{F:rtree}
\end{figure}

Note that the feasible sets for the deletion $G-a$ are $\{\emptyset,\{b\},\{c\},\{b,c\}\}$ and the feasible sets for the contraction $G/a$ are $\{\emptyset,\{b\},\{c\},\{d\}\}$.  In both cases, these are the independent sets in a matroid -- see the bottom of Figure~\ref{F:rtree} for geometric depictions of these two matroids.  Thus, it is possible for a greedoid $G$ that is not a matroid to have $G/x$ and $G-x$ both be matroids.

Note that $r(a)=1$ and $r(G-a)=r(G)=2$, so the recursion $$T(G)=(x-1)^{r(G)-r(G-a)}f(G-a)+(y-1)^{1-r(a)}f(G/a)$$ from Theorem~\ref{T:tutte-gen} simplifies to the familiar matroid recursion  $T(G)=T(G-a)+T(G/a)$.  Now $T(G-a)=x^2y$ since $G-a$ consists of 2 isthmuses and 1 loop, while $T(G/a)=x+y+y^2$ since $G/a\cong U_{1,3}$.  Thus  $T(G;x,y)=x^2 y+x+y+y^2.$

In this case, note that there is no matroid $M$ with $M-p\cong G-a$ and $M/p \cong G-a$ for some $p \in S$.  What's more, it is straightforward to check that there is no matroid with this Tutte polynomial. This gives an example of a greedoid that is not a matroid whose Tutte polynomial has  all positive coefficients.   
\end{exmp}

Although there is no matroid $M$ satisfying $T(M)=T(G)$ for the greedoid $G$ from Example~\ref{Ex:pos}, we can also construct examples where $T(G)=T(M)$ for a greedoid $G$ that is not a matroid and a matroid $M$.

\begin{exmp}\label{Ex:same tutte}
 Let $M$ be the matroid with geometric representation shown on the left in Figure~\ref{F:k4}, so $r(M)=3$ and $|S|=6$.  Let $G$ be the greedoid on the ground set $S=\{a,b,c,d,e,f\}$ with feasible sets $\emptyset$, all singletons, all pairs \textit{except } $ab$ and all triples  \textit{except }  $abc, ade, bef, cdf$.  One way to depict the collection of feasible triples is the set of all subsets of three labeled edges of $K_4$, except for the four subsets incident to a vertex, as in the graph on the right in Figure~\ref{F:k4}.  Then it is routine to verify $G$ is a greedoid.

\begin{figure}[h]
\centerline{\includegraphics[width=3in]{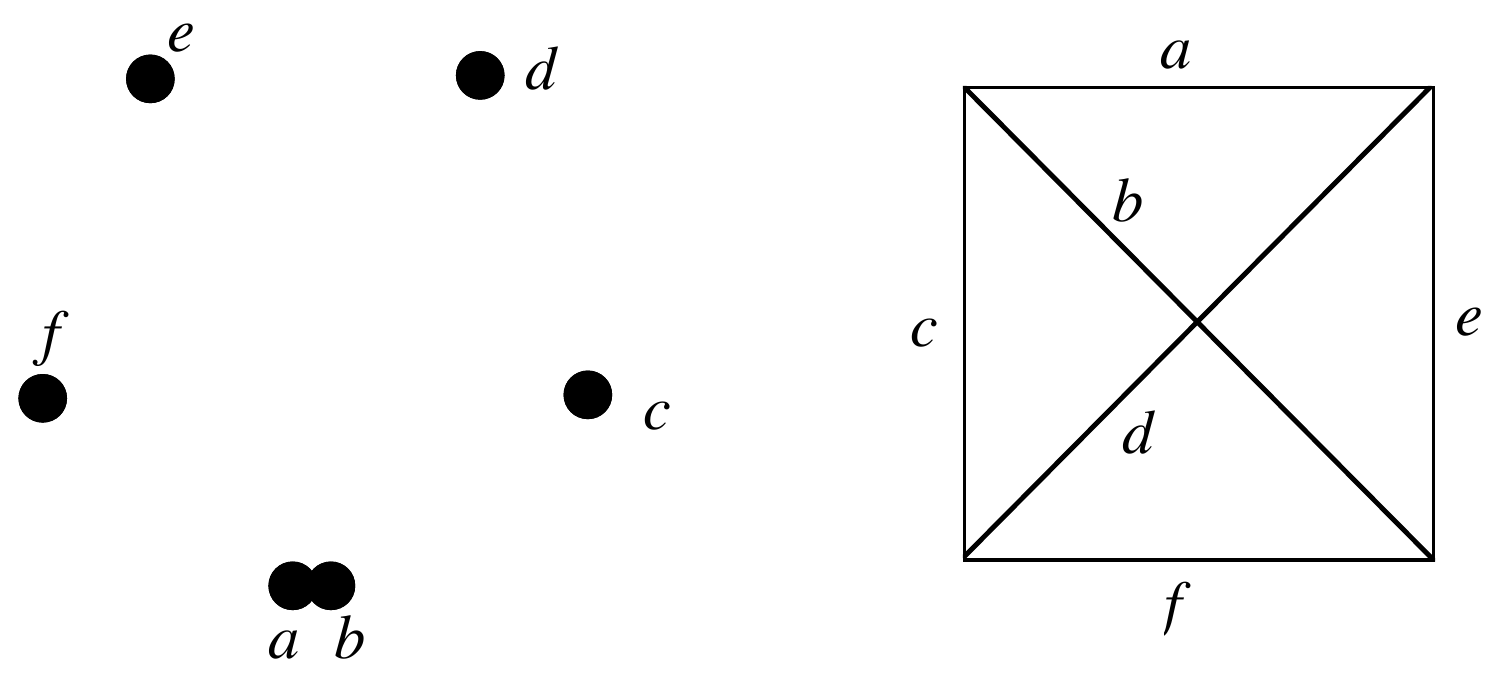}}
\caption{Left: Matroid $M$. Right: Edge labeled $K_4$ used to describe the feasible sets of the greedoid $G$. See Example~\ref{Ex:same tutte}.}
\label{F:k4}
\end{figure}

Now a computation shows
$$
T(M)=T(G)= x^3+x^2 y+2 x^2+2 x y+3 x+y^3+3 y^2+3 y
$$

In this case, there is no $p$ with both  $G-p$ and $G/p$ matroids.  Thus, if $T(G_1)=T(G_2)$ for greedoids $G_1$ and $G_2$, it need not be true that $G_1-a_1 \cong G_2-a_2$ and $G_1/a_1\cong G_2/a_2$ for some $a_1 \in S_1$ and $a_2\in S_2$.
\end{exmp}

In Example~\ref{Ex:pos}, we saw that $G-a$ and $G/a$ were both matroids, but $G-a$ contains a greedoid loop.  In that example, there is no matroid $M$ with $T(G)=T(M)$.  In Example~\ref{Ex:same tutte}, we found a matroid $M$ and a greedoid $G$ (where $G$ is not a matroid) with the same Tutte polynomial. Examples of two matroids (or greedoids) with the same Tutte polynomial abound, but virtually all of these arise from instances where $G_1-a_1 \cong G_2-a_2$ and $G_1/a_1\cong G_2/a_2$ for some $a_1 \in S_1$ and $a_2\in S_2$.

The next example gives such an example for a matroid $M$ and a greedoid $G$.  Thus, it combines features of Examples~\ref{Ex:pos} and \ref{Ex:same tutte}.

\begin{exmp}\label{Ex:same tutte2}
 Let $M'=M-f$ be the matroid obtained by deleting $f$ from the matroid $M$ of Example~\ref{Ex:same tutte}, so $r(M)=3$ and $|S|=5$.  Let $G$ be the greedoid on the ground set $S=\{a,b,c,d,e\}$ with feasible sets $\emptyset$, all singletons, all pairs \textit{except } $ab$ and all triples  \textit{except }  $abc, abd$ and $cde$. 
 
 Then $M'-e=G-e$ and $M'/e=G/e$, so $T(M')=T(G)$.  Hence, even when $G/p$ and $G-p$ are both matroids, with $r(G-p)=r(G)$ where $G-p$ has no loops, $G$ need not be a matroid.
 \end{exmp}

We conclude with an example that examines the role of duality for greedoids. Given a ranked set $G$, we can use Def.~\ref{D:gens}(1) to define a dual object $G^*$ from the rank function. When $G$ is a matroid, this agrees with the usual definition of duality. But matroid duals can also be defined  from the bases of $G^*$; these are simply the complements of the bases of $G$.

 While these two approaches are identical for matroids, they are not equivalent for greedoids. Using the rank formulation for duality, Theorem 4.4 of \cite{gdual} implies $G^*$ is also a greedoid iff $G$ is a  matroid. Thus, greedoid duality from the rank function does not exist for greedoids that are not also matroids. But it is frequently the case that the basis complements of a  greedoid $G$ form the bases for a different, ``dual-like'' greedoid (often in several distinct ways). 

But even this weak form of duality does not hold for all greedoids. This is the point of our final example.

\begin{exmp}\label{Ex:820} Exercise 8.20 of \cite{bz} asks for an example of a greedoid whose basis complements do not form the bases of any greedoid.  Let $G$ be the truncation of the branching greedoid associated with the rooted graph of Fig.~\ref{F:820}. Then $G$ has feasible sets 
$$\mathcal{F}=\{\emptyset, \{a\}, \{b\}, \{a,b\},\{a,c\},\{b,d\},\{a,b,c\},\{a,b,d\},\{a,c,e\},\{b,d,e\}\}.$$
\begin{figure}[h]
\centerline{\includegraphics[width=1.1in]{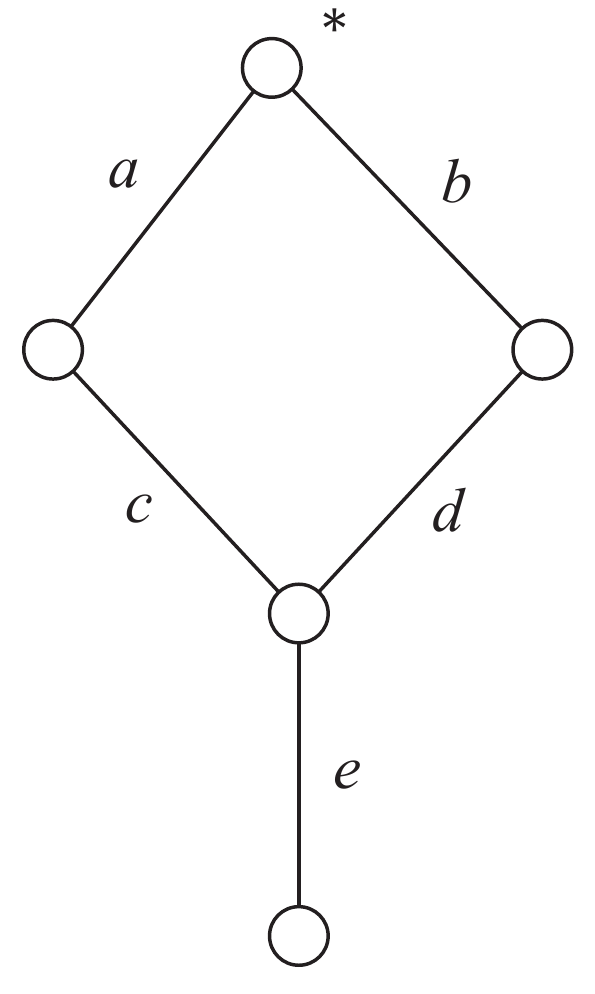}}
\caption{$G$ is the  truncation of the branching greedoid associated with the rooted tree above. The collection of basis complements cannot be the bases of any greedoid.}
\label{F:820}
\end{figure}

Then the set of bases complements of $G$ is $\{\{a,c\},\{b,d\},\{c,e\},\{d,e\}\}.$ To see that there is no greedoid having this collection as its bases, note that $\{a\}$ cannot be a feasible set (else feasible augmentation fails for the pair of feasible sets $\{a\}$ and $\{b,d\}$). Similarly, $\{c\}$ cannot be a feasible set. But then the feasible set $\{a,c\}$ is inaccessible.

We can compute the Tutte polynomial of $G$ (and recall $T(G^*;x,y)=T(G;y,x)$ from Theorem~\ref{T:tutte-gen}(2)).
$$T(G;x,y)=x^3 y^3-3 x^2 y^3+2 x^2 y^2+3 x y^3-4 x y^2+3 x y-y^3+3 y^2.$$

Although greedoids cannot be characterized by their bases -- for instance, an antimatroid has a unique basis -- it is frequently the case that there is a matroid whose bases coincide with the bases of the greedoid. But this example also demonstrates that the set of bases of a greedoid do not, in general,  satisfy the matroid basis properties (if they did, then, by matroid duality, their complements would, also).
\end{exmp}

\end{document}